\documentclass[a4paper,oneside]{article}%
\usepackage{CJK,amsmath}
\usepackage{amsfonts}
\usepackage{amsthm}
\usepackage{indentfirst}
\usepackage{slashbox}
\usepackage{makeidx}
\usepackage{mathrsfs}
\usepackage{amsmath}
\usepackage{amssymb}
\usepackage{amsfonts}
\usepackage{undertilde}
\usepackage[perpage,symbol*]{footmisc}
\usepackage{bm}
\usepackage{graphicx}
\usepackage[T1]{fontenc}
\usepackage{titlesec}
\usepackage{ccaption}
\usepackage[tight,hang]{subfigure}
\usepackage{graphicx}%
\providecommand{\U}[1]{\protect\rule{.1in}{.1in}}
\input{undertilde}
\makeatletter

\newcommand{\Rmnum}[1]{\expandafter\@slowromancap\romannumeral #1@}
\makeatother
\newcommand*{\HEI}
{\CJKfamily{hei}}

\newtheorem{theorem}{Theorem}[section]
\newtheorem{lemma}{Lemma}[section]

\newtheorem{corollary}{Corollary}[section]
\newtheorem{definition}{Definition}[section]

\newtheorem{assumption}{Assumption}[section]
\begin{document}

\title{Classical Solutions of Path-dependent PDEs and Functional Forward-Backward
Stochastic Systems\thanks{This work was supported by National Natural Science
Foundation of China (No. 11171187, No. 10871118 and No. 10921101).}}
\author{Shaolin Ji\thanks{Institute for Financial Studies and Institute of
Mathematics, Shandong University, Jinan, Shandong 250100, PR China
(Jsl@sdu.edu.cn, Fax: +86 0531 88564100).}
\and Shuzhen Yang\thanks{School of mathematics, Shandong University, Jinan,
Shandong 250100, PR China. (yangsz@mail.sdu.edu.cn). }}
\date{}
\maketitle

\textbf{Abstract}. In this paper we study the relationship between functional
forward-backward stochastic systems and {path-dependent PDEs. In the framework
of functional It\^{o} calculus, we introduce a path-dependent PDE and prove
that its solution is uniquely determined by a functional }forward-backward
stochastic system.

\bigskip

{\textbf{Keywords}: Functional It\^{o} calculus, Functional forword-backward
systems, Path-dependent PDEs, Classical solutions}.

\addcontentsline{toc}{section}{\hspace*{1.8em}Abstract}

\section{Introduction}

It is well known that quasilinear parabolic partial differential equations are
related to Markovian forward-backward stochastic differential equations (see
\cite{Peng S 5}, \cite{Pardoux.E 2} and \cite{Pardoux.E}), which generalizes
the classical Feynman-Kac formula. Recently in the framework of {functional
It\^{o} calculus, }a path-dependent PDE was introduced by Dupire
\cite{Dupire.B} and the so-called functional Feynman-Kac formula was also
obtained. For a recent account and development of this theory we refer the
reader to \cite{Cont.R}, \cite{Cont-2}, \cite{Cont-3}, \cite{Peng S 3},
{\cite{Peng S 2}, \cite{E}} and \cite{Cosso}.

In this paper, we study a {functional forword-backward system and its relation
to a quasilinear parabolic }path-dependent PDE. In more details, the{
functional forword-backward system is described by the following
forword-backward} SDE:
\begin{equation}
X^{\gamma_{t},x}(s)=\gamma_{t}(t)+\int_{t}^{s}b(X_{r}^{\gamma_{t}})dr+\int%
_{t}^{s}\sigma(X_{r}^{\gamma_{t}})dW(r), \tag{1.1}%
\end{equation}%
\begin{equation}
Y^{\gamma_{t}}(s)=g(X_{T}^{\gamma_{t}})-\int_{s}^{T}h(X_{r}^{\gamma_{t}%
},Y^{\gamma_{t}}(r),Z^{\gamma_{t}}(r))dr-\int_{s}^{T}Z^{\gamma_{t}%
}(r)dW(r),\quad s\in\lbrack t,T]. \tag{1.2}%
\end{equation}
After establishing some estimates and regularity results for the solution with
respect to paths, we prove that the solution of (1.2) is the unique classical
solution of the following path-dependent PDE
\begin{align*}
&  D_{t}u(\gamma_{t})+\mathcal{L}u(\gamma_{t})=h(\gamma_{t},u(\gamma
_{t}),D_{x}u(\gamma_{t})\sigma(\gamma_{t})),\\
&  u(\gamma_{T})=g(\gamma_{T}),\quad\gamma_{T}\in{\Lambda}.
\end{align*}
where
\[
\mathcal{L}u=\frac{1}{2}tr[\sigma\sigma^{T}D_{xx}u]+\langle b,D_{x}u\rangle.
\]

The paper is organized as follows: in section 2, we give the notations and
results on functional It\^{o} calculus and functional SDEs. Some estimates and
regularity results for the solution of FBSDEs are established in section 3.
Finally, we prove the relationship between functional FBSDEs and
{path-dependent PDEs} in section 4.

\section{Preliminaries}

\subsection{Functional It\^{o} calculus}

The following notations and tools are mainly from Dupire \cite{Dupire.B}. Let
$T>0$ be fixed. For each $t\in\lbrack0,T]$, we denote by $\Lambda_{t}$ the set
of c\`{a}dl\`{a}g $\mathbb{R}^{d}$-valued functions on $[0,t]$. For each
$\gamma\in\Lambda_{T}$ the value of $\gamma$ at time $s\in\lbrack0,T]$ is
denoted by $\gamma(s)$. Thus $\gamma=\gamma(s)_{0\leq s\leq T}$ is a
c\`{a}dl\`{a}g process on $[0,T]$ and its value at time $s$ is $\gamma(s)$.
The path of $\gamma$ up to time $t$ is denoted by $\gamma_{t}$, i.e.,
$\gamma_{t}=\gamma(s)_{0\leq s\leq t}\in\Lambda_{t}$. We denote $\Lambda
=\bigcup_{t\in\lbrack0,T]}\Lambda_{t}$. For each $\gamma_{t}\in\Lambda$ and
$x\in\mathbb{R}^{d}$ we denote by $\gamma_{t}(s)$ the value of $\gamma_{t}$ at
$s\in\lbrack0,t]$ and $\gamma_{t}^{x}:=(\gamma_{t}(s)_{0\leq s<t},\gamma
_{t}(t)+x)$ which is also an element in $\Lambda_{t}$.

Let $(\cdot,\cdot)$ and $|\cdot|$ denote the inner product and norm in
$\mathbb{R}^{n}$. We now define a distance on $\Lambda$. For each $0\leq
t,\bar{t}\leq T$ and $\gamma_{t},\bar{\gamma}_{\bar{t}}\in\Lambda$, we denote
\begin{align*}
\Vert\gamma_{t}\Vert:  &  =\sup\limits_{s\in\lbrack0,t]}|\gamma_{t}(s)|,\\
\Vert\gamma_{t}-\bar{\gamma}_{\bar{t}}\Vert:  &  =\sup\limits_{s\in
\lbrack0,t\vee\bar{t}]}|\gamma_{t}(s\wedge t)-\bar{\gamma}_{\bar{t}}%
(s\wedge\bar{t})|,\\
d_{\infty}(\gamma_{t},\bar{\gamma}_{\bar{t}}):  &  =\sup_{0\leq s\leq
t\vee\bar{t}}|\gamma_{t}(s\wedge t)-\bar{\gamma}_{\bar{t}}(s\wedge\bar
{t})|+|t-\bar{t}|.
\end{align*}
It is obvious that $\Lambda_{t}$ is a Banach space with respect to $\Vert
\cdot\Vert$ and $d_{\infty}$ is not a norm.

\begin{definition}
A function $u:\Lambda\mapsto\mathbb{R}$ is said to be $\Lambda$--continuous at
$\gamma_{t}\in\Lambda$, if for any $\varepsilon>0$ there exists $\delta>0$
such that for each $\bar{\gamma}_{\bar{t}}\in\Lambda$ with $d_{\infty}%
(\gamma_{t},\bar{\gamma}_{\bar{t}})<\delta$, we have $|u(\gamma_{t}%
)-u(\bar{\gamma}_{\bar{t}})|<\varepsilon$. $u$ is said to be $\Lambda
$--continuous if it is $\Lambda$--continuous at each $\gamma_{t}\in\Lambda$.
\end{definition}

\begin{definition}
Let $u:\Lambda\mapsto\mathbb{R}$ and $\gamma_{t}\in\Lambda$ be given. If there
exists $p\in\mathbb{R}^{d}$, such that
\[
u(\gamma_{t}^{x})=u(\gamma_{t})+\langle p,x\rangle+o(|x|)\ \text{as}%
\ x\rightarrow0, \ x\in\mathbb{R}^{d}.\ \
\]
Then we say that $u$ is (vertically) differentiable at $\gamma_{t}$ and denote
the gradient of $D_{x}u(\gamma_{t})=p$. $u$ is said to be vertically
differentiable in $\Lambda$ if $D_{x}u(\gamma_{t})$ exists for each
$\gamma_{t}\in\Lambda$. We can similarly define the Hessian $D_{xx}%
u(\gamma_{t})$. It is an $\mathbb{S}(d)$-valued function defined on $\Lambda$,
where $\mathbb{S}(d)$ is the space of all $d\times d$ symmetric matrices.
\end{definition}

For each $\gamma_{t}\in\Lambda$ we denote
\[
\gamma_{t,s}(r)=\gamma_{t}(r)\mathbf{1}_{[0,t)}(r)+\gamma_{t}(t)\mathbf{1}%
_{[t,s]}(r),\ \ r\in\lbrack0,s].
\]
It is clear that $\gamma_{t,s}\in\Lambda_{s}$.

\begin{definition}
For a given $\gamma_{t}\in\Lambda$ if we have
\[
u(\gamma_{t,s})=u(\gamma_{t})+a(s-t)+o(|s-t|)\ \text{as}\ s\rightarrow t,
\ s\geq t,\ \
\]
then we say that $u(\gamma_{t})$ is (horizontally) differentiable in $t$ at
$\gamma_{t}$ and denote $D_{t}u(\gamma_{t})=a$. $u$ is said to be horizontally
differentiable in $\Lambda$ if $D_{t}u(\gamma_{t})$ exists for each
$\gamma_{t}\in\Lambda$.
\end{definition}

\begin{definition}
Define $\mathbb{C}^{j,k}(\Lambda)$ as the set of function $u:=(u(\gamma
_{t}))_{\gamma_{t}\in\Lambda}$ defined on $\Lambda$ which are $j$ times
horizontally and $k$ times vertically differentiable in $\Lambda$ such that
all these derivatives are $\Lambda$--continuous.
\end{definition}

The following It\^{o} formula was firstly obtained by Dupire \cite{Dupire.B}
and then generalized by Cont and Fourni\'{e} \cite{Cont.R}, \cite{Cont-2} and
\cite{Cont-3}.

\begin{theorem}
[Functional It\^{o}'s formula]\label{w2} Let $(\Omega,\mathcal{F}%
,(\mathcal{F}_{t})_{t\in\lbrack0,T]},P)$ be a probability space, if $X$ is a
continuous semi-martingale and $u$ is in $\mathbb{C}^{1,2}(\Lambda)$, then for
any $t\in\lbrack0,T)$,
\[%
\begin{split}
u(X_{t})-u(X_{0})=  &  \int_{0}^{t}D_{s}u(X_{s})\,ds+\int_{0}^{t}D_{x}%
u(X_{s})\,dX(s)\\
&  +\frac{1}{2}\int_{0}^{t}D_{xx}u(X_{s})\,d\langle X\rangle(s)\quad
\quad\ P-a.s.
\end{split}
\]

\end{theorem}

\subsection{Functional FBSDEs}

Let $\Omega=C([0,T];\mathbb{R}^{d})$ and $P$ the Wiener measure on
$(\Omega,\mathbb{B}(\Omega))$. We denote by $W=(W(t)_{t\in\lbrack0,T]})$ the
cannonical Wiener process, with $W(t,\omega)=\omega(t)$, $t\in\lbrack0,T]$,
$\omega\in\Omega$. For any $t\in\lbrack0,T]$. we denote by $\mathcal{F}_{t}$
the $P$-completion of $\sigma(W(s),s\in\lbrack0,t])$.

For any $t\in\lbrack0,T]$, we denote by $L^{2}(\Omega,\mathcal{F}%
_{t};\mathbb{R}^{n})$ the set of all square integrable $\mathcal{F}_{t}%
-$measurable random variables, $M^{2}(0,T;\mathbb{R}^{n})$ the set of all
$\mathbb{R}^{n}$-valued $\mathcal{F}_{t}$-adapted processes $\vartheta(\cdot)$
such that
\[
E\int_{0}^{T}\mid\vartheta(s)\mid^{2}ds<+\infty.
\]

Let $t\in\lbrack0,T]$ and $\gamma_{t}\in{\Lambda}$. For every $s\in\lbrack
t,T]$, we consider the following functional forward-backward SDEs:
\begin{equation}
X^{\gamma_{t}}(s)=\gamma_{t}(t)+\int_{t}^{s}b(X_{r}^{\gamma_{t}})dr+\int%
_{t}^{s}\sigma(X_{r}^{\gamma_{t}})dW(r), \tag{2.1}%
\end{equation}%
\begin{equation}
Y^{\gamma_{t}}(s)=g(X_{T}^{\gamma_{t}})-\int_{s}^{T}h(X_{r}^{\gamma_{t}%
},Y^{\gamma_{t}}(r),Z^{\gamma_{t}}(r))dr-\int_{s}^{T}Z^{\gamma_{t}}(r)dW(r),
\tag{2.2}%
\end{equation}
where%
\[
X^{\gamma_{t}}(s):=\gamma_{t}(s),\,s\in\lbrack0,t].
\]
The processes $X,Y,Z$ take values in $\mathbb{R}^{n},\mathbb{R}^{n}%
,\mathbb{R}^{n\times d}$; $b,h,\sigma$ and $g$ take values in $\mathbb{R}%
^{n},\mathbb{R}^{n},\mathbb{R}^{n\times d}$ and $\mathbb{R}^{n}$. $(2.1)$ and
$(2.2)$ can be rewriten as:
\begin{align*}
&  dX^{\gamma_{t}}(s)=b(X_{s}^{\gamma_{t}})ds+\sigma(X_{s}^{\gamma_{t}%
})dW(s),\\
&  dY^{\gamma_{t}}(s)=h(X_{s}^{\gamma_{t}},Y^{\gamma_{t}}(s),Z^{\gamma_{t}%
}(s))ds+Z^{\gamma_{t}}(s)dW(s),\\
&  X^{\gamma_{t}}(t)=\gamma_{t}(t),\qquad Y^{\gamma_{t}}(T)=g(X_{T}%
^{\gamma_{t}}).
\end{align*}

For $z\in\mathbb{R}^{n\times d}$, we define $|z|=\{tr(zz^{T})\}^{1/2}$.
\noindent For $z^{1}\in\mathbb{R}^{n\times d}$, $z^{2}\in\mathbb{R}^{n\times
d}$,
\[
((z^{1},z^{2}))=tr(z^{1}(z^{2})^{T}),
\]
and for $u^{1}=(y^{1},z^{1})\in\mathbb{R}^{n}\times\mathbb{R}^{n\times d}$,
$u^{2}=(y^{2},z^{2})\in\mathbb{R}^{n}\times\mathbb{R}^{n\times d}$
\[
\lbrack u^{1},u^{2}]=\langle y^{1},y^{2}\rangle+((z^{1},z^{2})).
\]
%for $u=(x,y,z)\in \mathbb{R}^n \times \mathbb{R}^n \times \mathbb{R}^{n\times d}$,
%\begin{displaymath}
%f(u)=(h(u),b(u),\sigma(u)).
%\end{displaymath}
We give the following assumption:

\begin{assumption}
$\forall x_{t}^{1},$ $x_{t}^{2}\in{\Lambda}$ and $t\in\lbrack0,T],$ there
exists a constant $c_{1}>0$, such that%
\[
\mid b(x_{t}^{1})-b(x_{t}^{2})\mid+\mid\sigma(x_{t}^{1})-\sigma(x_{t}^{2}%
)\mid\leq c_{1}\parallel x_{t}^{1}-x_{t}^{2}\parallel,\quad a.e.
\]
and $\forall x_{t}\in{\Lambda,}$
\[
\mid b(x_{t})\mid+\mid\sigma(x_{t})\mid\leq c_{1}(1+\parallel x_{t}%
\parallel),\quad a.e.
\]

\end{assumption}

\begin{definition}
$X:[0,T]\times\Omega$$\rightarrow$ $\mathbb{R}^{n}$ is called an adapted
solution of the Eqs. $(2.1)$, if $X\in M^{2}(0,T;\mathbb{R}^{n})$, and it
satisfies $(2.1)$ $P-a.s.$
\end{definition}

Then we have the following theorem (see \cite{Lipster}):

\begin{theorem}
Let Assumptions $2.1$ hold, then there exists a unique adapted solution $X$
for Eqs. $(2.1)$.
\end{theorem}

\section{Regularity}

We first recall some notions in Pardoux and Peng \cite{Pardoux.E 2}.
$C^{n}(\mathbb{R}^{p};\mathbb{R}^{q})$, $C_{b}^{n}(\mathbb{R}^{p}%
;\mathbb{R}^{q})$, $C_{p}^{n}(\mathbb{R}^{p};\mathbb{R}^{q})$ will denote
respectively the set of functions of class $C^{n}$ from $\mathbb{R}^{p}$ into
$\mathbb{R}^{q}$, the set of those functions of class $C_{b}^{n}$ whose
partial derivatives of order less than or equal to $n$ are bounded, and the
set of those functions of class $C_{p}^{n}$ which, together with all their
partial derivatives of order less than or equal to n, grow at most like a
polynomial function of the variable $x$ at infinity.

Now we give the definition of derivatives in our context. Under the above
Assumption 2.1 we have
\begin{align*}
&  dX^{\gamma_{t}}(s)=b(X_{s}^{\gamma_{t}})ds+\sigma(X_{s}^{\gamma_{t}%
})dW(s),\\
&  X^{\gamma_{t}}(t)=\gamma_{t}(t),
\end{align*}
has a uniqueness solution. For $t\leq s\leq T$, set
\begin{align*}
\tilde{\Lambda}_{\gamma_{t},s}  &  :=\{\bar{\gamma}_{s}:\bar{\gamma
}(h)=X^{\gamma_{t}}(h,\omega),\;0\leq h\leq s,\,\omega\in\Omega\},\\
\tilde{\Lambda}_{t,s}  &  :=\bigcup_{\gamma_{t}\in\Lambda_{t}}\tilde{\Lambda
}_{\gamma_{t},s}\text{ and }\tilde{\Lambda}_{t}:=\bigcup_{0\leq s\leq T}%
\tilde{\Lambda}_{t,s}.
\end{align*}
Then the following definition of derivatives will be used frequently in the sequel.

\begin{definition}
An $\mathbb{R}^{n}$-valued function $g$ is said to be in $C^{2}(\tilde
{\Lambda}_{\gamma_{t},T})$, if for $\gamma_{1}\in\tilde{\Lambda}_{\gamma
_{t},T}$ and $\gamma_{2}\in\tilde{\Lambda}_{\gamma_{t}^{y},T}$, there exist
$p_{1}\in\mathbb{R}^{d}$ and $p_{2}\in\mathbb{R}^{d}\times\mathbb{R}^{d}$ such
that $p_{2}$ is symmetric,
\[
g(\gamma_{2})-g(\gamma_{1})=\langle p_{1},y\rangle+\frac{1}{2}\langle
p_{2}y,y\rangle+o(|y|^{2}),\quad x\in\mathbb{R}^{d}.\quad
\]
We denote $g{^{\prime}}_{\gamma_{t}}(\gamma_{1}):=p_{1}$, and $g{^{\prime
\prime}}_{\gamma_{t}}(\gamma_{1}):=p_{2}$. $g$ is said to be in $C_{l,lip}%
^{2}(\tilde{\Lambda}_{t,T})$ if $g{^{\prime}}_{\gamma_{t}}(\gamma)$ and
$g{^{\prime\prime}}_{\gamma_{t}}(\gamma)$ exist for each $\gamma_{t}%
\in{\Lambda_{t}}$, and there exists some constants $C\geq0$ and $k\geq0$
depending only on $g$ such that for each $\gamma,\bar{\gamma}\in\Lambda
_{T},t,s\in\lbrack0,T]$,
\[
\mid g(\gamma)-g(\bar{\gamma})\mid\leq C(\parallel\gamma\parallel
^{k}+\parallel\bar{\gamma}\parallel^{k})\parallel\gamma-\bar{\gamma}%
\parallel,
\]
and for each $\gamma\in\tilde{\Lambda}_{t,T},\bar{\gamma}\in\tilde{\Lambda
}_{s,T},t,s\in\lbrack0,T]$,
\[
\mid\Phi_{\gamma_{t}}(\gamma)-\Phi_{\gamma_{s}}(\bar{\gamma})\mid\leq
C(\parallel\gamma\parallel^{k}+\parallel\bar{\gamma}\parallel^{k})(\mid
t-s\mid+\parallel\gamma-\bar{\gamma}\parallel)
\]
with $\Phi=g{^{\prime}}_{\gamma_{t}}(\gamma),g{^{\prime\prime}}_{\gamma_{t}%
}(\gamma)$. We can also define $C^{2}(\tilde{\Lambda}_{t,s})$, $C_{l,lip}%
^{2}(\tilde{\Lambda}_{t,s})$, $C_{l,lip}^{1}(\tilde{\Lambda}_{t,s})$,
$C_{l,lip}(\tilde{\Lambda}_{t,s})$ and $C^{2}(\tilde{\Lambda}_{t})$,
$C_{l,lip}^{2}(\tilde{\Lambda}_{t})$, $C_{l,lip}^{1}(\tilde{\Lambda}_{t})$,
$C_{l,lip}(\tilde{\Lambda}_{t})$.
\end{definition}

Now we consider the solvability of equation (2.2).

\begin{assumption}
Let $g$ is an $\mathbb{R}^{n}$-valued function on $\Lambda_{T}$. Moreover
$g\in C_{l,lip}^{2}(\tilde{\Lambda}_{t,T})$ with the Lipschitz constants $C$
and $k$.
\end{assumption}

\begin{assumption}
Let $h(\gamma_{t},y,z)=\bar{h}(t,\gamma(t),y,z)$, where $\bar{h}%
:[0,T]\times\mathbb{R}^{n}\times\mathbb{R}^{n}\times\mathbb{R}^{n\times
d}\mapsto\mathbb{R}^{n}$ is such that $(t,r,y,z)\mapsto\bar{\Psi}(t,r,y,z)$ is
of class $C_{p}^{0,3}([0,T]\times\mathbb{R}^{n}\times\mathbb{R}^{n}%
\times{\mathbb{R}^{n\times d}};\mathbb{R}^{n})$ and the first order partial
derivatives in r,y and z are bounded, as well as their derivatives of up to
order two with respect to y,z.
\end{assumption}

It is obvious under Assumption 2.1, 3.1, 3.1 the FBSDE (2.1) and (2.2) has a
uniqueness solution (see \cite{Lipster}, \cite{Pardoux.E 3} and \cite{Ma J-1}).

\subsection{Regularity of the solution of FBSDEs}

We assume the Lipschitz constants with respect to $b,\sigma,h$ are C and k.
Then we have the following estimates for the solution of FBSDE (2.1) and
(2.2).\newline

\begin{lemma}
Under Assumption $2.1$, $3.1$, $3.2$ there exists $C_{2}$ and $q$ depending
only on $C,T,k,x$ such that
\begin{align*}
&  E[\sup_{s\in\lbrack t,T]}\mid X^{\gamma_{t}}(s)\mid^{2}]\leq C_{2}%
(1+\parallel\gamma_{t}\parallel^{2}),\\
&  E[\sup_{s\in\lbrack t,T]}\mid Y^{\gamma_{t}}(s)\mid^{2}]\leq C_{2}%
(1+\parallel\gamma_{t}\parallel^{q}),\\
&  E[(\int_{t}^{T}\mid Z^{\gamma_{t}}(s)\mid^{2}ds)]\leq C_{2}(1+\parallel
\gamma_{t}\parallel^{q}).
\end{align*}

\end{lemma}

\noindent\textbf{{Proof.}} To simplify presentation, we only study the case
$n=d=1$.\newline Applying It\^{o}'s formula to $(Y_{\gamma_{t},x}%
(s))^{2}e^{\beta_{1}s}$ yields that%
\[%
\begin{array}
[c]{rl}
& (Y^{\gamma_{t},}(s))^{2}e^{\beta_{1}s}+\int_{s}^{T}e^{\beta_{1}r}%
[(Z^{\gamma_{t}}(r))^{2}+\beta_{1}(Y^{\gamma_{t}}(r))^{2}]dr\\
= & g^{2}(X_{T}^{\gamma_{t}})e^{\beta_{1}T}-\int_{s}^{T}2e^{\beta_{1}%
r}Y^{\gamma_{t}}(r)h(X_{r}^{\gamma_{t}},Y^{\gamma_{t}}(r),Z^{\gamma_{t}%
}(r))dr-\int_{s}^{T}2e^{\beta_{1}r}Y^{\gamma_{t}}(r)Z^{\gamma_{t}}(r)dW(r).
\end{array}
\]
So%
\[%
\begin{array}
[c]{rl}
& (Y^{\gamma_{t}}(s))^{2}+E[\int_{s}^{T}e^{\beta_{1}(r-s)}[(Z^{\gamma_{t}%
}(r))^{2}+\beta_{1}(Y^{\gamma_{t}}(r))^{2}]dr\mid\mathcal{F}_{s}]\\
= & E[g^{2}(X_{T}^{\gamma_{t}})e^{\beta_{1}(T-s)}\mid\mathcal{F}_{s}%
]-E[\int_{s}^{T}2e^{\beta_{1}(r-s)}Y^{\gamma_{t}}(r)h(X_{r}^{\gamma_{t}%
},Y^{\gamma_{t}}(r),Z^{\gamma_{t}}(r))dr\mid\mathcal{F}_{s}].
\end{array}
\]
Then we have%
\[%
\begin{array}
[c]{rl}
& E\sup_{t\leq s\leq T}(Y^{\gamma_{t}}(s))^{2}+E[\int_{t}^{T}e^{\beta
_{1}(r-t)}[(Z^{\gamma_{t}}(r))^{2}+\beta_{1}(Y^{\gamma_{t}}(T))^{2}]dr]\\
\leq & E[g^{2}(X_{T}^{\gamma_{t}})e^{\beta_{1}(T-t)}]+E[\int_{t}^{T}%
e^{\beta_{1}(r-t)}\frac{2}{\beta_{1}}h^{2}(X_{r}^{\gamma_{t}},Y^{\gamma_{t}%
}(r),Z^{\gamma_{t}}(r))dr]+E[\int_{t}^{T}e^{\beta_{1}(r-t)}\frac{\beta_{1}}%
{2}(Y^{\gamma_{t}}(r))^{2}d(r)].
\end{array}
\]
and%
\[%
\begin{array}
[c]{rl}
& E\sup_{t\leq s\leq T}(Y^{\gamma_{t}}(s))^{2}+E[\int_{t}^{T}e^{\beta
\tag{3.3}_{1}(r-t)}[(Z^{\gamma_{t}}(r))^{2}+\frac{\beta_{1}}{2}(Y^{\gamma_{t}%
}(r))^{2}]dr]\\
\leq & E[g^{2}(X_{T}^{\gamma_{t}})e^{\beta_{1}(T-t)}]+E[\int_{t}^{T}%
e^{\beta_{1}(r-t)}\frac{2}{\beta_{1}}h^{2}(X_{r}^{\gamma_{t}},Y^{\gamma_{t}%
}(r),Z^{\gamma_{t}}(r))dr].
\end{array}
\]
Applying It\^{o}'s formula to $(X^{\gamma_{t},x}(s))^{2}$ yields that%
\[
(X^{\gamma_{t}}(s))^{2}={\gamma_{t}(t)}^{2}+\int_{t}^{s}2X^{\gamma_{t}%
}(r)b(X_{r}^{\gamma_{t}})dr+\int_{t}^{s}2X^{\gamma_{t}}(r)\sigma(X_{r}%
^{\gamma_{t}})dW(r)+\int_{t}^{s}\sigma^{2}(X_{r}^{\gamma_{t}})dr
\]
By inequality $2ab\leq a^{2}+b^{2}$ and Burkholder-Davis-Gundy's inequality,
there is a $C_{0}$ such that,
\[
E\sup_{t\leq r\leq s}(X^{\gamma_{t}}(r))^{2}\leq C_{0}[{\gamma_{t}(t)}%
^{2}+E\int_{t}^{s}b^{2}(X_{r}^{\gamma_{t}})dr+E\int_{t}^{s}(X^{\gamma_{t}%
}(r))^{2}dr+E\int_{t}^{s}\sigma^{2}(X_{r}^{\gamma_{t}})dr].
\]
By Assumption $2.1$ and Gronwall's inequality, we have (note that $C_{0}$ will
change line by line)
\[
E\sup_{t\leq r\leq T}(X^{\gamma_{t}}(r))^{2}\leq C_{0}(1+\parallel\gamma
_{t}\parallel^{2}).
\]
By Assumptions $3.1$ and $3.2$ and taking $\beta_{1}=4C^{2}+1$, we have%
\[%
\begin{array}
[c]{rl}
& E[\sup_{t\leq s\leq T}(Y^{\gamma_{t}}(s))^{2}+E[\int_{t}^{T}[(Z^{\gamma_{t}%
}(r))^{2}+(Y^{\gamma_{t}}(r))^{2}]dr]\\
\leq & C_{0}(1+\parallel\gamma_{t}\parallel^{q})
\end{array}
\]
where $q=2(1+k)$. This completes the proof.\quad$\Box$

Now we study the regularity properties of the solution of FBSDE (2.1), (2.2)
with respect to the "parameter" $\gamma_{t}$. For $0\leq s<t\leq T$, define
$Y^{\gamma_{t}}(s)=Y^{\gamma_{t}}(s\vee t)$ and $Z^{\gamma_{t}}(s)=0$.

\begin{theorem}
Under Assumptions $2.1$, $3.2$ and $3.3$, there exist $C_{2}$ and $q$
depending only on $C,c_{2},x$ such that for any $t,\bar{t}\in\lbrack0,T]$,
$\gamma_{t},\bar{\gamma}_{\bar{t}}$, and $h,\bar{h}\in{\mathbb{R}%
\setminus\{0\}}$.
\end{theorem}

\noindent$(i)$
\[
E[\sup_{u\in\lbrack t\vee\bar{t},T]}\mid Y^{\gamma_{t}}(u)-Y^{\bar{\gamma
}_{\bar{t}}}(u)\mid^{2}]\leq C_{2}(1+\parallel\gamma_{t}\parallel
^{q}+\parallel\bar{\gamma}_{\bar{t}}\parallel^{q})(\parallel\gamma_{t}%
-\bar{\gamma}_{\bar{t}}\parallel^{2}+\mid t-\bar{t}\mid),\newline%
\]
$(ii)$
\[
E[\int_{t\vee\bar{t}}^{T}\mid Z^{\gamma_{t}}(u)-Z^{\bar{\gamma}_{\bar{t}}%
}(u)\mid^{2}du]\leq C_{2}(1+\parallel\gamma_{t}\parallel^{q}+\parallel
\bar{\gamma}_{\bar{t}}\parallel^{q})(\parallel\gamma_{t}-\bar{\gamma}_{\bar
{t}}\parallel^{2}+\mid t-\bar{t}\mid),\newline%
\]
$(iii)$
\[%
\begin{array}
[c]{rl}
& E[\sup_{u\in\lbrack t\vee\bar{t},T]}\mid\Delta_{h}^{i}Y^{\gamma_{t}%
}(u)-\Delta_{h}^{i}Y^{\bar{\gamma}_{{\bar{t}}}}(u)\mid^{2}]\\
\leq & C_{2}(1+\parallel\gamma_{t}\parallel^{q}+\parallel\bar{\gamma}_{\bar
{t}}\parallel^{q}+\mid h\mid^{q}+\mid\bar{h}\mid^{q}))(\mid h-\bar{h}\mid
^{2}+\parallel\gamma_{t}-\bar{\gamma}_{\bar{t}}\parallel^{2}+\mid t-\bar
{t}\mid),\newline%
\end{array}
\]
$(iv)$
\[%
\begin{array}
[c]{rl}
& E[\int_{t\vee\bar{t}}^{T}\mid\Delta_{h}^{i}Z^{\gamma_{t}}(u)-\Delta_{h}%
^{i}Z^{\bar{\gamma}_{{\bar{t}}}}(u)\mid^{2}du]\\
\leq & C_{2}(1+\parallel\gamma_{t}\parallel^{q}+\parallel\bar{\gamma}_{\bar
{t}}\parallel^{q}+\mid h\mid^{q}+\mid\bar{h}\mid^{q}))(\mid h-\bar{h}\mid
^{2}+\parallel\gamma_{t}-\bar{\gamma}_{\bar{t}}\parallel^{2}+\mid t-\bar
{t}\mid),\newline%
\end{array}
\]
where
\[
\Delta_{h}^{i}Y^{\gamma_{t},x}(s)=\frac{1}{h}(Y^{\gamma_{t}^{h_{e_{i}}}%
}(s)-Y^{\gamma_{t}}(s)),\Delta_{h}^{i}Z^{\gamma_{t}}(s)=\frac{1}{h}%
(Z^{\gamma_{t}^{h_{e_{i}}}}(s)-Z^{\gamma_{t}}(s))
\]
and $(e_{1},\cdots,e_{n})$ is an orthonormal basis of $\mathbb{R}^{n}$.

\noindent\textbf{{Proof.}} $(Y^{\gamma_{t}}-Y^{\bar{\gamma}_{\bar{t}}%
},Z^{\gamma_{t}}-Z^{\bar{\gamma}_{\bar{t}}})$ can be formed as a linearized
BSDE: for each $s\in\lbrack t\vee\bar{t},T]$,%
\[%
\begin{array}
[c]{rl}
& Y^{\gamma_{t}}(s)-Y^{\bar{\gamma}_{\bar{t}}}(s)\\
= & g(X_{T}^{\gamma_{t}})-g(X_{T}^{\bar{\gamma}_{\bar{t}}})+\int_{s}%
^{T}[h(X_{r}^{\gamma_{t}},Y^{\gamma_{t}}(r),Z^{\gamma_{t}}(r))-h(X_{r}%
^{\bar{\gamma}_{\bar{t}}},Y^{\bar{\gamma}_{\bar{t}}}(r),Z^{\bar{\gamma}%
_{\bar{t}}}(r))]dr+\int_{s}^{T}(Z^{\gamma_{t}}(r)-Z^{\bar{\gamma}_{\bar{t}}%
}(r))dW(r)\\
= & g(X_{T}^{\gamma_{t}})-g(X_{T}^{\bar{\gamma}_{\bar{t}}},)-\int_{s}^{T}%
[\hat{\alpha}_{\gamma_{t},\bar{\gamma}_{\bar{t}}}(r)+\hat{\beta}_{\gamma
_{t},\bar{\gamma}_{\bar{t}}}(Y^{\gamma_{t}}(r)-Y^{\bar{\gamma}_{\bar{t}}%
}(r))+\hat{\delta}_{\gamma_{t},\bar{\gamma}_{\bar{t}}}(Z^{\gamma_{t}%
}(r)-Z^{\bar{\gamma}_{\bar{t}}}(r))]dr\\
& +\int_{s}^{T}(Z^{\gamma_{t}}(r)-Z^{\bar{\gamma}_{\bar{t}}}(r))]dW(r),
\end{array}
\]
where (with $U^{\gamma_{t}}=(Y^{\gamma_{t}},Z^{\gamma_{t}})$)%
\[%
\begin{array}
[c]{rl}%
\hat{\alpha}_{\gamma_{t},\bar{\gamma}_{\bar{t}}}(r) & =h(X_{r}^{{\gamma}_{{t}%
}},Y^{\bar{\gamma}_{\bar{t}}}(r),Z^{\bar{\gamma}_{\bar{t}}}(r))-h(X_{r}%
^{\bar{\gamma}_{\bar{t}}},Y^{\bar{\gamma}_{\bar{t}}}(r),Z^{\bar{\gamma}%
_{\bar{t}}}(r)),\\
\hat{\beta}_{\gamma_{t},\bar{\gamma}_{\bar{t}}}(r)(Y^{\gamma_{t}}%
(r)-Y^{\bar{\gamma}_{\bar{t}}}(r)) & =\int_{0}^{1}\frac{\partial h}{\partial
y}(X_{r}^{\gamma_{t}},U^{\bar{\gamma}_{\bar{t}}}(r)+\theta(U^{{\gamma}_{{t}}%
}(r)-U^{\bar{\gamma}_{\bar{t}}}(r)))d\theta,\\
\hat{\delta}_{\gamma_{t},\bar{\gamma}_{\bar{t}}}(r)(Z^{\gamma_{t}}%
(r)-Z^{\bar{\gamma}_{\bar{t}}}(r)) & =\int_{0}^{1}\frac{\partial h}{\partial
z}(X_{r}^{\gamma_{t}},U^{\bar{\gamma}_{\bar{t}}}(r)+\theta(U^{{\gamma}_{{t}}%
}(r)-U^{\bar{\gamma}_{\bar{t}}}(r)))d\theta.
\end{array}
\]
Under Assumptions $3.1$ , $3.2$, using the same method as in Lemma $3.1,$ we
get the first three inequalities.

For the next three inequalities, we write $(\Delta_{h}^{i}Y^{\gamma_{t}%
},\Delta_{h}^{i}Z^{\gamma_{t}})$ as the solution of the following linearized
BSDE:%
\[%
\begin{array}
[c]{rl}
& \Delta_{h}^{i}Y^{\gamma_{t}}(s)\\
= & \frac{1}{h}(g(X_{T}^{\gamma_{t}^{h_{e_{i}}}})-g(X_{T}^{\gamma_{{t}}%
}))-\int_{s}^{T}[\frac{1}{h}\hat{\alpha}_{{\gamma_{t},{\gamma}_{{t}}%
}^{h_{e_{i}}}}(r)+\hat{\beta}_{\gamma_{t},{{\gamma}_{{t}}}^{h_{e_{i}}}%
}(r)\Delta_{h}^{i}Y^{\gamma_{t}}(r)+\hat{\delta}_{\gamma_{t},{{\gamma}_{{t}}%
}^{h_{e_{i}}}}\Delta_{h}^{i}Z^{\gamma_{t}}(r)]dr\\
& -\int_{t}^{T}\Delta_{h}^{i}Z^{\gamma_{t}}(r)dW(r).
\end{array}
\]
Then the same calculus implies that
\[
E[\sup_{s\in\lbrack t,T]}\mid\Delta_{h}^{i}Y^{\gamma_{t}}(s)\mid^{2}+\mid
\int_{t}^{T}\mid\Delta_{h}^{i}Z^{\gamma_{t}}\mid^{2}dr\mid]\leq C_{2}%
(1+\parallel\gamma_{t}\parallel^{q}+\mid h\mid^{q}).
\]
Consider%
\[%
\begin{array}
[c]{rl}
& \Delta_{h}^{i}Y^{\gamma_{t}}(s)-\Delta_{\bar{h}}^{i}Y^{\bar{\gamma}_{\bar
{t}}}(s)\\
= & \frac{1}{h}(g(X_{T}^{\gamma_{t}^{h_{e_{i}}}})-g(X_{T}^{\gamma_{{t}}%
}))-\frac{1}{\bar{h}}(g(X_{T}^{\bar{\gamma}_{t}^{{\bar{h}}_{e_{i}}}}%
)-g(X_{T}^{\bar{\gamma}_{\bar{t}}}))-\int_{s}^{T}(\Delta_{h}^{i}Z^{\gamma_{t}%
}(r)-\Delta_{\bar{h}}^{i}Z^{\bar{\gamma}_{\bar{t}}}(r))dW(r)\\
& -\{\int_{s}^{T}[\frac{1}{h}\hat{\alpha}_{\gamma_{t},{\gamma}_{{t}}%
^{{{h}_{e_{i}}}}}(r)-\frac{1}{\bar{h}}\hat{\alpha}_{\bar{\gamma}_{\bar{t}%
},\bar{\gamma}_{\bar{t}}^{{\bar{h}_{e_{i}}}}}(r)+\hat{\beta}_{\gamma
_{t},{{\gamma}_{{t}}}^{{{h}_{e_{i}}}}}(r)\Delta_{{h}}^{i}Y^{\gamma_{t}%
}(r)-\hat{\beta}_{\bar{\gamma}_{\bar{t}},{\bar{\gamma}_{\bar{t}}^{{\bar
{h}_{e_{i}}}}}}(r)\Delta_{\bar{h}}^{i}Y^{\bar{\gamma}_{\bar{t}}}(r)\\
& +\hat{\delta}_{\gamma_{t},{\gamma}_{{t}}}\Delta_{{h}}^{i}Z^{\gamma_{t}%
}(r)-\hat{\delta}_{\bar{\gamma}_{\bar{t}},{\bar{\gamma}_{\bar{t}}^{{\bar
{h}_{e_{i}}}}}}(r)\Delta_{\bar{h}}^{i}Z^{\bar{\gamma}_{\bar{t}}}(r)]dr\}.
\end{array}
\]
Set
\[
(\tilde{Y}(s),\tilde{Z}(s)):=(\Delta_{h}^{i}Y^{\gamma_{t}}(s)-\Delta_{\bar{h}%
}^{i}Y^{\bar{\gamma}_{\bar{t}}}(s),\Delta_{h}^{i}Z^{\gamma_{t}}(s)-\Delta
_{\bar{h}}^{i}Z^{\bar{\gamma}_{\bar{t}}}(s)).
\]
Then it solves the following BSDE%
\[%
\begin{array}
[c]{rl}%
\tilde{Y}(s)= & \frac{1}{h}(g(X_{T}^{\gamma_{t}^{h_{e_{i}}}})-g(X_{T}%
^{\gamma_{{t}}}))-\frac{1}{\bar{h}}(g(X_{T}^{\bar{\gamma}_{t}^{{\bar{h}%
}_{e_{i}}}})-g(X_{T}^{\bar{\gamma}_{\bar{t}}}))\\
& -\int_{s}^{T}[{\hat{\beta}}_{\gamma_{t},{{\gamma}_{{t}}}^{{{h}_{e_{i}}}}%
}(r)\tilde{Y}(r)+{\hat{\delta}}_{\gamma_{t},{{\gamma}_{{t}}}^{{{h}_{e_{i}}}}%
}\tilde{Z}(r)+\tilde{h}(r)]dr-\int_{s}^{T}\tilde{Z}(r)dW(r),
\end{array}
\]
where
\[
\tilde{h}(r):=[{\hat{\beta}}_{\gamma_{t},{{\gamma}_{{t}}}^{{{h}_{e_{i}}}}%
}(r)-{\hat{\beta}}_{\bar{\gamma}_{\bar{t}},{\bar{\gamma}_{\bar{t}}^{{\bar
{h}_{e_{i}}}}}}(r)]\Delta_{\bar{h}}^{i}Y^{\bar{\gamma}_{\bar{t}}}%
(r)+[{\hat{\delta}}_{\gamma_{t},{{\gamma}_{{t}}}^{{{h}_{e_{i}}}}}%
(r)-{\hat{\delta}}_{\bar{\gamma}_{\bar{t}},{\bar{\gamma}_{\bar{t}}^{{\bar
{h}_{e_{i}}}}}}(r)]\Delta_{\bar{h}}^{i}Z^{\bar{\gamma}_{\bar{t}}}(r)+\frac
{1}{h}{\hat{\alpha}}_{\gamma_{t},{\gamma}_{{t}}^{{{h}_{e_{i}}}}}(r)-\frac
{1}{\bar{h}}{\hat{\alpha}}_{\bar{\gamma}_{\bar{t}},\bar{\gamma}_{\bar{t}%
}^{{\bar{h}_{e_{i}}}}}(r).
\]
Thus, under Assumptions 3.1, 3.2, similarly as in Lemma 3.1, we can get the
last three inequalities.\quad$\Box$

\begin{theorem}
For each $\gamma_{t}\in\Lambda$, $\{Y^{{\gamma_{t}^{z}}}(s),\,s\in\lbrack
t,T],z\in\mathbb{R}^{n}\}$ has a version which is a.e. of class $C^{0,2}%
([0,T]\times\mathbb{R}^{n})$.
\end{theorem}

\noindent\textbf{{Proof.}} We only consider one dimensional case. Applying
Lemma $3.1$, for each $h,\bar{h}\in\mathbb{R}\setminus\{0\}$ and $k,\bar{k}%
\in\mathbb{R}$,\newline%
\[%
\begin{array}
[c]{rl}
& E[\sup_{u\in\lbrack t,T]}\mid Y^{\gamma_{t}^{k}}(u)-Y^{{\gamma}_{{t}}%
^{\bar{k}}}(u)\mid^{2}]\leq C_{p}(1+\parallel\gamma_{t}\parallel^{q})\mid
k-\bar{k}\mid^{2},\\
& E[\mid\int_{t}^{T}\mid Z^{\gamma_{t}^{k}}(u)-Z^{{\gamma}_{{t}}^{\bar{k}}%
}\mid^{2}du\mid]\leq C_{2}(1+\parallel\gamma_{t}\parallel^{q})\mid k-\bar
{k}\mid^{2},\\
& E[\sup_{u\in\lbrack t,T]}\mid\Delta_{h}^{i}Y^{\gamma_{t}^{k}}(u)-\Delta
_{\bar{h}}^{i}Y^{{\gamma}_{{t}}^{\bar{k}}}\mid^{2}]\\
\leq & C_{2}(1+\parallel\gamma_{t}\parallel^{q}+\parallel\bar{\gamma}_{\bar
{t}}\parallel^{q}+\mid h\mid^{q}+\mid\bar{h}\mid^{q}))(\mid k-\bar{k}\mid
^{2}+\mid h-\bar{h}\mid^{2}),\\
& E[\mid\int_{t}^{T}\mid\Delta_{h}^{i}Z^{\gamma_{t}^{k}}(u)-\Delta_{\bar{h}%
}^{i}Z^{\gamma_{{t}}^{\bar{k}}}(u)\mid^{2}du\mid]\\
\leq & C_{2}(1+\parallel\gamma_{t}\parallel^{q}+\parallel\bar{\gamma}_{\bar
{t}}\parallel^{q}+\mid h\mid^{q}+\mid\bar{h}\mid^{q}))(\mid k-\bar{k}\mid
^{2}+\mid h-\bar{h}\mid^{2}).
\end{array}
\]
By kolmogorov's criterion, there exists a continuous derivative of
$Y^{\gamma_{t}^{z}}(s)$ with respect to $z.$ There also exists a mean-square
derivative of $Z^{\gamma_{t}^{z}}(s)$ with respect to $z$, which is mean
square continuous in $z$. We denote them by
\[
(D_{z}Y^{\gamma_{t}},D_{z}Z^{\gamma_{t}}).
\]
By Theorem 3.1 and definition $3.1$, $(D_{z}Y^{\gamma_{t}},D_{z}Z^{\gamma_{t}%
})$ is the solution of the following BSDE:%
\[%
\begin{array}
[c]{rl}%
D_{x}Y^{\gamma_{t}}(s)= & g_{\gamma_{t}}^{\prime}(X_{T}^{\gamma_{t}})-\int%
_{s}^{T}[h_{\gamma_{t}}^{\prime}(X_{r}^{\gamma_{t}},Y^{\gamma_{t}%
}(r),Z^{\gamma_{t}}(r))-\int_{s}^{T}D_{x}Z^{\gamma_{t}}(r)dW(r)\\
& +h_{y}^{\prime}(X_{r}^{\gamma_{t}},Y^{\gamma_{t}}(r),Z^{\gamma_{t}}%
(r))D_{x}Y^{\gamma_{t}}(r)+h_{z}^{\prime}(X_{r}^{\gamma_{t}},Y^{\gamma_{t}%
}(r),Z^{\gamma_{t}}(r))D_{x}Z^{\gamma_{t}}(r)]dr
\end{array}
\]
It is easy to check that the above BSDE has a uniqueness solutions. Thus the
existence of a continuous second order derivative of $Y^{\gamma_{t}^{z}}(s)$
with respect to $z$ is proved in a similar way. $\Box$\newline Define
\[
u(\gamma_{t}):=Y^{\gamma_{t}}(t),\quad for\quad\gamma_{t}\in\Lambda.
\]
We have the following results about $u(\gamma_{t})$.\quad

\begin{lemma}
$\forall t\leq s \leq T$, we have $u(X_{s}^{\gamma_{t}})=Y^{\gamma_{t}}(s)$.
\end{lemma}

\noindent\textbf{{Proof.}} For given $\gamma_{t_{1}}$, $t_{1}<t$, set
$X(r)=xI_{0\leq r\leq t_{1}}$. Consider the solution of FBSDE (2.1) and (2.2)
on $[t,T]$:%
\[%
\begin{tabular}
[c]{ll}%
$X^{\gamma_{t}}(s)=$ & $\gamma_{t}(t)+\int_{t}^{s}b(X_{r}^{\gamma_{t}}%
)dr+\int_{t}^{s}\sigma(X_{r}^{\gamma_{t}})dW(r),$\\
$Y^{\gamma_{t}}(s)=$ & $g(X_{T}^{\gamma_{t}})-\int_{s}^{T}h(X_{r}^{\gamma_{t}%
},Y^{\gamma_{t}}(r),Z^{\gamma_{t}}(r))dr-\int_{s}^{T}Z^{\gamma_{t}%
}(r)dW(r),\quad s\in\lbrack t,T].$%
\end{tabular}
\ \ \
\]
We need to prove $u(X_{t}^{\gamma_{t_{1}}})=Y^{\gamma_{t_{1}}}(t)$. Define%
\[
X_{t}^{N,\gamma_{t_{1}}}:=\Sigma_{i=1}^{N}I_{A_{i}}x_{t}^{i},
\]
where $\{A_{i}\}_{i=1}^{N}$ is a division of $\mathcal{F}_{t}$, $x_{t}^{i}\in
A_{i}\cap{\Lambda}$, $i=1,2,\cdots,N$. For any $i$, $(Y^{x_{t}^{i},a^{i}%
}(s),Y^{x_{t}^{i},a^{i}}(s))$ is the solution of the following BSDE:
\[
Y^{x_{t}^{i}}(s)=g(X_{T}^{x_{t}^{i}})-\int_{s}^{T}h(X_{r}^{x_{t}^{i}}%
,Y^{x_{t}^{i}}(r),Z^{x_{t}^{i}}(r))dr-\int_{s}^{T}Z^{x_{t}^{i}}(r)dW(r),\quad
s\in\lbrack t,T].
\]
Multiplying by $I_{A_{i}}$ and adding the corresponding terms, we obtain:%
\[%
\begin{array}
[c]{rl}%
\sum_{i=1}^{N}I_{A_{i}}Y^{x_{t}^{i}}(s)= & g(\sum_{i=1}^{N}I_{A_{i}}%
X_{T}^{x_{t}^{i}})-\int_{s}^{T}h(\sum_{i=1}^{N}I_{A_{i}}X_{r}^{x_{t}^{i}}%
,\sum_{i=1}^{N}I_{A_{i}}Y^{x_{t}^{i}}(r),\sum_{i=1}^{N}I_{A_{i}}Z^{x_{t}^{i}%
}(r))dr\\
& -\int_{s}^{T}\sum_{i=1}^{N}I_{A_{i}}Z^{x_{t}^{i}}(r)dW(r),\quad s\in\lbrack
t,T].
\end{array}
\]
By the uniqueness and existence theorem of BSDE, we get $Y_{s}^{\gamma_{t}%
}=\sum_{i=1}^{N}I_{A_{i}}Y^{x_{t}^{i}}(s)$, $Z_{s}^{\gamma_{t}}=\sum_{i=1}%
^{N}I_{A_{i}}Z^{x_{t}^{i}}(s)$ a.s. Then, by the definition of $u$, we get
\[
Y^{{\gamma_{t_{1}}}}(t)=\sum_{i=1}^{N}I_{A_{i}}Y^{x_{t}^{i}}(t)=\sum_{i=1}%
^{N}I_{A_{i}}u(x_{t}^{i})=u(X_{t}^{\gamma_{t_{1}}}).
\]
For the general case, following the method in Peng and Wang \cite{Peng S 3}
(Lemma $4.3$), we choose a simple adapted process $\{\gamma_{t}^{i}%
\}_{i=1}^{\infty}$ such that $E\parallel\gamma_{t}^{i}-X_{t}^{\gamma_{t_{1}}%
}\parallel$ convergence to $0$ as $i\rightarrow\infty$. We obtain\newline%
\[
E\mid Y^{\gamma_{t}^{i}}(t)-Y^{{\gamma_{t_{1}}}}(t)\mid^{2}\leq CE\parallel
\gamma_{t}^{i}-X_{t}^{\gamma_{t_{1}}}\parallel
\]
This completes the proof.\quad$\Box$

By Theorem 3.1 and 3.2 and the definition of vertical derivative, we have the
following corollary.

\begin{corollary}
$u(\gamma_{t})$ is $\Lambda$-continuous and $D_{x}u(\gamma_{t})$%
,$D_{xx}u(\gamma_{t}) $ exist, moreover they are both $\Lambda$-continuous.
\end{corollary}

\noindent\textbf{{Proof.}} By Theorem $3.2$ we know that $D_{z}u(\gamma_{t})$
and $D_{zz}u(\gamma_{t})$ exist. In the following, we only prove $u(\gamma
_{t})$ is $\Lambda$-continuous. The proof for the continuous property of
$D_{z}u(\gamma_{t})$ and $D_{zz}u(\gamma_{t})$ is similar. Taking expectation
on both sides of equation (2.2),
\begin{equation}
u({\gamma_{t}})=Eg(X_{T}^{\gamma_{t}})-E\int_{t}^{T}h(X_{r}^{\gamma_{t}%
},Y^{\gamma_{t}}(r),Z^{\gamma_{t}}(r))dr. \tag{5.9}%
\end{equation}
For $\gamma_{t},\bar{\gamma}_{\bar{t}}\in\Lambda,$ $\bar{t}\geq t$, we have%
\[%
\begin{array}
[c]{rl}
& \mid u(\gamma_{t})-u(\bar{\gamma}_{\bar{t}})\mid\\
\leq & E[\mid g(X_{T}^{\gamma_{t}})-g(X_{T}^{\bar{\gamma}_{\bar{t}}}%
)\mid]+E[\int_{t}^{\bar{t}}\mid h(X_{r}^{\gamma_{t}},Y^{\gamma_{t}%
}(r),Z^{\gamma_{t}}(r))\mid dr]\\
& +E[\int_{\bar{t}}^{T}\mid h(X_{r}^{\gamma_{t}},Y^{\gamma_{t}}(r),Z^{\gamma
_{t}}(r))-h(X_{r}^{\bar{\gamma}_{\bar{t}}},Y^{\bar{\gamma}_{\bar{t}}%
}(r),Z^{\bar{\gamma}_{\bar{t}}}(r))\mid dr]\\
\leq & E[C_{1}(1+\parallel X_{T}^{\gamma_{t}}\parallel^{k}+\parallel
X_{T}^{\bar{\gamma}_{\bar{t}}}\parallel^{k})\parallel{\gamma_{t}}-{\bar
{\gamma}_{\bar{t}}}\parallel\\
& +3(\bar{t}-t)^{\frac{1}{2}}(\int_{t}^{\bar{t}}(\mid h(X^{\bar{\gamma}%
_{\bar{t}}}(r),0,0)\mid^{2}+\mid CY^{\bar{\gamma}_{\bar{t}}}(r)\mid^{2}\\
& +\mid CY^{\bar{\gamma}_{\bar{t}}}(r)\mid^{2})dr)^{\frac{1}{2}}+C\int%
_{\bar{t}}^{T}(\mid Y^{\gamma_{t}}(r)-Y^{\bar{\gamma}_{\bar{t}}}\mid+\mid
Z^{\gamma_{t}}(r)-Z^{\bar{\gamma}_{\bar{t}}}\mid)dr].
\end{array}
\]
By Theorem 3.1, for some constant $C_{1}$ depending only in $C,k$ and $T$,
\[
\mid u(\gamma_{t})-u(\bar{\gamma}_{\bar{t}})\mid\leq C_{1}(1+\parallel
\gamma_{t}\parallel^{k}+\parallel\bar{\gamma}_{\bar{t}}\parallel
^{k})(\parallel\gamma_{t}-\bar{\gamma}_{\bar{t}}\parallel+\mid t-\bar{t}%
\mid^{\frac{1}{2}}).
\]
This completes the proof.\quad$\Box$

\subsection{Path regularity of process Z}

In Pardoux and Peng \cite{Pardoux.E 2}, BSDE is only state-dependent, i.e.,
$h=h(t,\gamma(t),y,z)$ and $g=g(\gamma(T))$. Under appropriate assumptions,
$Y$ and $Z$ are related in the following sense:
\[
Z^{\gamma_{t}}(s)=\nabla_{x}u(s,\gamma_{t}(t)+W(s)-W(t)),\quad P-a.s.
\]
Peng and Wang \cite{Peng S 3} extends this result to the path-dependent case.
The corresponding BSDE is
\[
Y^{\gamma_{t}}(s)=g(W_{T}^{\gamma_{t}})-\int_{s}^{T}h(W_{r}^{\gamma_{t}%
},Y^{\gamma_{t}}(r),Z^{\gamma_{t}}(r))dr-\int_{s}^{T}Z^{\gamma_{t}%
}(r)dW(r),\quad s\in\lbrack t,T].
\]
where $W_{T}^{\gamma_{t}}=I_{s\leq t}\gamma_{t}(s)+I_{t<s\leq T}(\gamma
_{t}(t)+W(s)-W(t))$. Then under some assumptions, they obtained
\[
Z^{\gamma_{t}}(s)=D_{x}u(W_{s}^{\gamma_{t}}),\quad P-a.s.
\]
In our context, we have the following theorem:

\begin{theorem}
Under Assumption $2.1$, $3.1$ and $3.2$, for each $\gamma_{t}\in\Lambda$, the
process $(Z^{\gamma_{t}}(s))_{s\in\lbrack t,T]}$ has a continuous version with
the form,
\[
Z^{\gamma_{t}}(s)=\sigma(X_{s}^{\gamma_{t}})D_{x}u(X_{s}^{\gamma_{t}}),\quad
for\quad s\in\lbrack t,T]\quad P-a.s.
\]

\end{theorem}

To prove the above Theorem, we need the following lemma essentially from
Pardoux and Peng \cite{Pardoux.E 2}.

\begin{lemma}
Let $\gamma_{t}$ and some $\bar{t}\in\lbrack t,T]$ be given. Suppose that
\[
g(\gamma,z)=\varphi(\gamma(\bar{t}),\gamma(T)-\gamma(\bar{t}),z),
\]
where $\varphi$ is in $C_{p}^{3}(\mathbb{R}^{2d}\times\mathbb{R}%
^{m};\mathbb{R}^{m})$. For $\phi=b,\sigma$ and $h$, suppose that
\[
h(\gamma_{t},y,z)=h_{1}(s,\gamma_{s}(s),y,z)I_{[0,\bar{t})}(s)+h_{2}%
(s,\gamma_{s}(s)-\gamma_{s}(\bar{t}),y,z)I_{[\bar{t},T]}(s),
\]
where $h_{1},h_{2}\in C^{0,3}$. Then for each $s\in\lbrack t,T]$,
\[
Z^{\gamma_{t}}(s)=\sigma(X_{s}^{\gamma_{t}})D_{x}u(X_{s}^{\gamma_{t}}),\quad
for\quad s\in\lbrack t,T]\quad\ P-a.s.
\]

\end{lemma}

\noindent\textbf{{Proof.}} We only consider the one dimensional case. For
$s\in\lbrack\bar{t},T]$, the BSDE $(2.2)$ can be rewritten as%
\[%
\begin{array}
[c]{rl}%
Y^{\gamma_{t}}(u)= & \varphi(\gamma_{s}(\bar{t}),X^{\gamma_{s}}({T}%
)-\gamma_{s}(\bar{t}))-\int_{u}^{T}Z^{\gamma_{s}}(r)dW(r)\\
& -\int_{u}^{T}{h}_{2}(r,\gamma_{s}(\bar{t}),X^{\gamma_{s}}(r)-\gamma_{s}%
(\bar{t}),Y^{\gamma_{s}}(r),Z^{\gamma_{s}}(r))dr,\quad u\in\lbrack s,T].
\end{array}
\]
For $s\in\lbrack t,\bar{t}]$,%
\[%
\begin{array}
[c]{rl}%
Y^{\gamma_{s}}(u)= & \varphi(X^{\gamma_{s}}(\bar{t}),X^{\gamma_{s}}%
({T})-X^{\gamma_{s}}(\bar{t}))-\int_{u}^{T}Z^{\gamma_{s}}(r)dW(r)\\
& -\int_{u}^{T}{h}_{2}(r,X^{\gamma_{s}}(\bar{t}),X^{\gamma_{s}}(r)-X^{\gamma
_{s}}(\bar{t}),Y^{\gamma_{s}}(r),Z^{\gamma_{s}}(r))dr,\quad u\in\lbrack\bar
{t},T],
\end{array}
\]%
\[
Y^{\gamma_{s}}(u)=Y^{\gamma_{s}}(\bar{t})-\int_{u}^{\bar{t}}{h}_{1}%
(r,X^{\gamma_{s}}(\bar{t}),Y^{\gamma_{s}}(r),Z^{\gamma_{s}}(r))dr-\int%
_{u}^{\bar{t}}Z^{\gamma_{s}}(r)dW(r),\quad u\in\lbrack s,\bar{t}].
\]
Now consider the following system of quasilinear parabolic differential
equations, which is defined on $[\bar{t},T]\times\mathbb{R}^{2}$ and
parameterized by $x\in\mathbb{R}$,
\begin{align*}
&  \partial_{s}u_{2}(s,x,y)+\mathcal{L}u_{2}(s,x,y)={h}_{2}(s,x,y,u_{2}%
(s,x,y),\partial_{y}u_{2}(s,x,y)\sigma(r_{s})\\
&  u_{2}(T,x,y)=\varphi(x,y).
\end{align*}
where $\mathcal{L}=\frac{1}{2}\sigma^{2}\frac{\partial^{2}}{\partial_{xx}%
}+b\frac{\partial}{\partial_{x}}$. The other one is defined on $[t,\bar
{t}]\times\mathbb{R}$:
\begin{align*}
&  \partial_{s}u_{1}(s,x)+\mathcal{L}u_{1}(s,x)={h}_{1}(r,x,u_{1}%
(s,x),\partial_{y}u_{1}(s,x)\sigma(r)\\
&  u_{1}(\bar{t},x)=u_{2}(\bar{t},x,0).
\end{align*}
where $\mathcal{L}=\frac{1}{2}\sigma^{2}\frac{\partial^{2}}{\partial_{xx}%
}+b\frac{\partial}{\partial_{x}}$. By the above Corollary 3.1 and Theorem
$3.1$ , $3.2$ of Paroux-Peng \cite{Pardoux.E 2}, we have $u_{2}\in
C^{1,2}([\bar{t},T]\times\mathbb{R}^{2};\mathbb{R})$, $u_{1}\in C^{1,2}%
([t,\bar{t}]\times\mathbb{R};\mathbb{R})$ and
\[
u(\gamma_{s})=u_{1}(s,\gamma_{s}(s))I_{[t,\bar{t})}(s)+u_{2}(s,\gamma_{s}%
(\bar{t}),\gamma_{s}(s)-\gamma_{s}(\bar{t}))I_{[\bar{t},T]}(s).
\]
Then we obtain
\begin{align*}
&  Y^{\gamma_{t}}(s)=u_{1}(s,X^{\gamma_{t}}(s)),\quad t\leq s<\bar{t},\\
&  Y^{\gamma_{t}}(s)=u_{2}(s,X^{\gamma_{t}}(\bar{t}),X^{\gamma_{t}%
}(s)-X^{\gamma_{t}}(\bar{t})),\quad\bar{t}\leq s\leq T,\\
&  Z^{\gamma_{t}}(s)=\partial_{x}u_{1}(s,X^{\gamma_{t}}(s))\sigma
(X_{s}^{\gamma_{t}}),\quad t\leq s<\bar{t},\\
&  Z^{\gamma_{t}}(s)=\partial_{x}u_{2}(s,X^{\gamma_{t}}(\bar{t}),X^{\gamma
_{t}}(s)-X^{\gamma_{t}}(\bar{t}))\sigma(X_{s}^{\gamma_{t}}),\quad\bar{t}\leq
s\leq T.
\end{align*}

Finally, for each $s\in\lbrack t,T]$,
\[
\sigma(X_{s}^{\gamma_{t}})D_{x}u(X_{s}^{\gamma_{t}})=Z^{\gamma_{t}}%
(s)\quad\ \ P-a.s.
\]
In particular,
\[
\sigma(\gamma_{t})D_{x}u(\gamma_{t})=Z^{\gamma_{t}}(t).\quad\gamma_{t}%
\in\Lambda.
\]
This completes the proof.\quad$\Box$

Now we give the proof of Theorem 3.3.

\noindent\textbf{{Proof.}} For each fixed $t\in\lbrack0,T]$ and positive
integer $n$, we introduce a mapping $\gamma^{n}(\bar{\gamma_{s}}):\Lambda
_{s}\mapsto\Lambda_{s}$
\[
\gamma^{n}(\bar{\gamma}_{s})(r)=\bar{\gamma}_{s}(r)I_{[0,t)}+\Sigma
_{k=0}^{n-1}\bar{\gamma}_{s}(t_{k+1}^{n}\wedge s)I_{[t_{k}^{n}\wedge
s]}(r)+\bar{\gamma}_{s}(s)I_{\{s\}}(r),\,s\in\lbrack0,T],
\]
where $t_{k}^{n}=t+\frac{k(T-t)}{n}$, $k=0,1,\ldots,n$
\[
g^{n}(\bar{\gamma}):=g(\gamma^{n}(\bar{\gamma})),\quad h^{n}(\bar{\gamma}%
_{s},y,z):=h(\gamma^{n}(\bar{\gamma_{s}}),y,z).
\]
For each $n$, there exists some functions $\varphi_{n}$ defined on
$\Lambda_{t}\times\mathbb{R}^{n\times d}$ and $\psi_{n}$ defined on
$[t,T]\times{\Lambda_{t}}\times\mathbb{R}^{n\times d}\times{\mathbb{R}%
^{m}\times\mathbb{R}^{m\times d}}$ such that%
\[%
\begin{array}
[c]{rl}%
g^{n}(\bar{\gamma})= & \varphi_{n}(\bar{\gamma}_{t},\bar{\gamma}(t_{1}%
^{n})-\bar{\gamma}(t),\cdots,\bar{\gamma}(t_{n}^{n})-\bar{\gamma}(t_{n-1}%
^{n})),\\
h^{n}(\bar{\gamma}_{s},y,z)= & \psi_{n}(s,\bar{\gamma}_{t},\bar{\gamma}%
_{s}(t_{1}^{n}\wedge s)-\bar{\gamma}_{s}(t),\cdots,\bar{\gamma}_{s}(t_{n}%
^{n}\wedge s)-\bar{\gamma}_{s}(t_{n-1}^{n}\wedge s),y,z).
\end{array}
\]
Indeed, if we set
\begin{align*}
&  \bar{\varphi}_{n}(\bar{\gamma}_{t},x_{1},\cdots,x_{n}):=g((\bar{\gamma}%
_{t}(s)I_{[0,t)}(s)+\Sigma_{k=1}^{n}x_{k}I_{[t_{k-1}^{n},t_{k}^{n})}%
(s)+x_{n}I_{\{T\}}(s))_{0\leq s\leq T}),\\
&  \varphi_{n}(\bar{\gamma}_{t},x_{1},\cdots,x_{n}):=\bar{\varphi}_{n}%
(\bar{\gamma}_{t},\bar{\gamma}_{t}+x_{1},\bar{\gamma}_{t}(t)+x_{1}%
+x_{2},\cdots,\bar{\gamma}_{t}(t)+\Sigma_{i=1}^{n}x_{i}),
\end{align*}
then by Assumptions 2.1, 3.1 and 3.2, we obtain that, for each fixed
$\bar{\gamma}_{t}$, $\varphi_{n}(\bar{\gamma}_{t},x_{1},\cdots,x_{n})$ is a
$C_{p}^{3}$-function of $x_{1},\cdots,x_{n}$. In particular, for each
$\bar{\gamma}\in\Lambda$,
\[
\partial_{x_{i}}\varphi_{n}(\bar{\gamma}_{t},\bar{\gamma}(t_{1}^{n}%
)-\bar{\gamma}(t),\cdots,\bar{\gamma}(t_{n}^{n})-\bar{\gamma}(t_{n-1}%
^{n}))=g_{\gamma_{t_{i-1}^{n}}}^{\prime}(\gamma^{n}(\bar{\gamma})).
\]
\newline

For any $\bar{t}\geq t$, $\bar{\gamma}_{\bar{t}}\in\Lambda_{\bar{t}}$, we
consider the following BSDE:
\[
Y^{n,\bar{\gamma}_{\bar{t}}}(s)=g^{n}(X_{T}^{n,\bar{\gamma}_{\bar{t}}}%
)-\int_{s}^{T}h^{n}(X_{r}^{n,\bar{\gamma}_{\bar{t}}},Y^{n,\bar{\gamma}%
_{\bar{t}}}(r),Z^{n,\bar{\gamma}_{\bar{t}}}(r))dr-\int_{s}^{T}Y^{n,\bar
{\gamma}_{\bar{t}}}(r)dW(r).
\]
we denote
\[
u^{n}(\bar{\gamma}_{\bar{t}}):=Y^{n,\bar{\gamma}_{\bar{t}}}(t),\quad
\bar{\gamma}_{\bar{t}}\in\Lambda.
\]
Following the argument as in Lemma 3.3, for each $s\in\lbrack t,T]$, we have
\[
\sigma^{n}(X_{s}^{n,\bar{\gamma}_{\bar{t}}})D_{x}u^{n}(\bar{\gamma}_{\bar{t}%
})=Z^{n,\bar{\gamma}_{\bar{t}}}(s)\quad a.s.
\]
Let $C_{0}$ be a constant depending only on $C,T$ and $k$, which is allowed to
change from line by line. Following the similar calculus as in Lemma $3.1$ and
Theorem $3.1,$ we get that%

\[%
\begin{array}
[c]{cl}
& u^{n}(\bar{\gamma}_{\bar{t}})-u(\bar{\gamma}_{\bar{t}})\\
= & \gamma^{n}(\bar{\gamma}_{\bar{t}})(\bar{t})-\bar{\gamma}_{\bar{t}}(\bar
{t})+g^{n}(X_{T}^{\bar{\gamma}_{\bar{t}}})-g(X_{T}^{\bar{\gamma}_{\bar{t}}})\\
& +\int_{\bar{t}}^{T}[h^{n}(X_{r}^{\bar{\gamma}_{\bar{t}}},Y^{n,\bar{\gamma
}_{\bar{t}}}(r),Z^{n,\bar{\gamma}_{\bar{t}}}(r))-h(X_{r}^{\bar{\gamma}%
_{\bar{t}}},Y^{\bar{\gamma}_{\bar{t}}}(r),Z^{\bar{\gamma}_{\bar{t}}%
}(r))]dr+\int_{\bar{t}}^{T}[Z^{n,\bar{\gamma}_{\bar{t}}}(r))-Z^{\bar{\gamma
}_{\bar{t}}}(r)]dW(r)
\end{array}
\]
\newline For $X_{s}^{n,\bar{\gamma}_{t}}=\gamma^{n}(X_{T}^{\bar{\gamma}_{t}}%
)$, we have the next result,%
\[%
\begin{array}
[c]{rl}%
\lim_{n}X_{T}^{n,\bar{\gamma}_{\bar{t}}} & =X_{T}^{\bar{\gamma}_{\bar{t}}%
},\text{ }P-a.s.\\
\lim_{n}(Y^{n,\bar{\gamma}_{\bar{t}}}(s),Z^{n,\bar{\gamma}_{\bar{t}}}(s)) &
=(Y^{\bar{\gamma}_{\bar{t}}}(s),Z^{\bar{\gamma}_{\bar{t}}}(s)),\text{
}a.e.s\in\lbrack t,T],\text{ }P-a.s.
\end{array}
\]
We can get%
\[
\lim_{n}u^{n}(\bar{\gamma}_{\bar{t}})=u(\bar{\gamma}_{\bar{t}}),\text{ }%
\lim_{n}{D_{x}}u^{n}(\bar{\gamma}_{\bar{t}})={D_{x}}u(\bar{\gamma}_{\bar{t}%
}),\text{ }\lim_{n}{D_{xx}}u^{n}(\bar{\gamma}_{\bar{t}})={D_{xx}}u(\bar
{\gamma}_{\bar{t}}),
\]
and \bigskip%
\[%
\begin{array}
[c]{rl}%
\lim_{n}(u^{n}(X_{s}^{n,\bar{\gamma}_{t}}),{D_{x}}u^{n}(X_{s}^{n,\bar{\gamma
}_{t}}),{D_{xx}}u^{n}(X_{s}^{n,\bar{\gamma}_{t}})) & =(u(X_{s}^{\bar{\gamma
}_{t}}),{D_{x}}u(X_{s}^{\bar{\gamma}_{t}}),{D_{xx}}u(X_{s}^{\bar{\gamma}_{t}%
})),
\end{array}
\]
\indent $\text{ }a.e.s\in\lbrack t,T],\text{ }P-a.s.$ \newline So that
\[
Z^{\gamma_{t}}(s)=\sigma(X_{s}^{\bar{\gamma}_{t}})D_{x}u(X_{s}^{\bar{\gamma
}_{t}}),\ \ \text{\ }a.e.s\text{ }\in\lbrack t,T],\text{ }P-a.s.
\]
\quad This completes the proof. $\Box$

\section{The related path-dependent PDEs}

In this section, we relate FBSDE (2.1), (2.2) to the following path-dependent
partial differential equation:
\begin{align*}
&  D_{t}u(\gamma_{t})+\mathcal{L}u(\gamma_{t})-h(\gamma_{t},u(\gamma
_{t}),\sigma(\gamma_{t})D_{x}u(\gamma_{t}))=0,\\
&  u(\gamma_{T})=g(\gamma_{T}),\quad\gamma_{T}\in{\Lambda}^{n}.
\end{align*}
where
\[
\mathcal{L}u=\frac{1}{2}tr[(\sigma\sigma^{T})D_{xx}u]+\langle b,D_{x_{{}}%
}u\rangle.
\]

\begin{theorem}
Suppose Assumption $2.1$, $3.1$ and $3.2$ hold, and if $u\in\mathbb{C}%
^{1.2}(\Lambda)$ and $u$ is the solutions of equation $(4.1)$, $u$ is
uniformly Lipschitz continuous, and bounded by $C(1+\parallel\gamma
_{t}\parallel)$, then the solution is uniqueness, and for any $\gamma_{t}%
\in\Lambda$, $u(\gamma_{t})$ is determined by equation $(2.1)$ and $(2.2)$.
\end{theorem}

\noindent\textbf{{Proof. }}By the assumptions of this theorem,, we know that
$b(\gamma_{t})$ and $\sigma(\gamma_{t})$ is uniformly Lipschtiz continuous and
the following SDE has a uniqueness solution.%
\[%
\begin{array}
[c]{rl}%
dX^{\gamma_{t}}(s)= & b(X_{s}^{\gamma_{t}})ds+\sigma(X_{s}^{\gamma_{t}%
})dW(s),\\
X_{t}= & \gamma_{t},\quad s\in\lbrack t,T].
\end{array}
\]
Set $Y(s)=u({X_{s}^{\gamma_{t}}}),\quad t\leq s\leq T$ , Applying It\^{o}'s
formula to $Y(s)=u({X_{s}^{\gamma_{t}}})$, we have%
\[%
\begin{array}
[c]{rl}%
dY(s)= & -h(X_{s}^{\gamma_{t}},Y(s),Z(s))dr-\sigma(X_{s}^{\gamma_{t}}%
)D_{x}u(X_{s}^{\gamma_{t}})dW(s),\\
Y(T)= & g(X_{T}^{\gamma_{t}})\quad s\in\lbrack t,T].
\end{array}
\]
Then by the uniqueness and existence theorem of the functional FBSDE, we
obtain the result.\quad$\Box$

Now we prove the converse to the about result.

\begin{theorem}
Under Assumption $2.1$, $3.1$ and $3.2$. The function $u(\gamma_{t}%
)=Y^{\gamma_{t}}(t)$ is the unique $\mathbb{C}^{1,2}(\Lambda)$-solution of the
path-dependent PDE $(4.1)$.
\end{theorem}

\noindent\textbf{{Proof.}} We only study the one dimensional case.
$u\in\mathbb{C}^{0,2}(\Lambda)$ follows from Corollary $3.1$. Let $\delta>0$
be such that $t+\delta\leq T$. By Lemma $3.2$ we can get
\[
u(X_{t+\delta}^{\gamma_{t}})=Y^{\gamma_{t}}(t+\delta).
\]
Hence
\[
u(\gamma_{t,t+\delta})-u(\gamma_{t})=u(\gamma_{t,t+\delta})-u(X_{t+\delta
}^{\gamma_{t}})+u(X_{t+\delta}^{\gamma_{t}})-u(\gamma_{t}),
\]
By the proof of Theorem $3.3$, we obtain
\[
u(\gamma_{t,t+\delta})-u(\gamma_{t})=\lim_{n\rightarrow\infty}[u^{n}%
(\gamma_{t,t+\delta})-u^{n}(X_{t+\delta}^{\gamma_{t}})]+\int_{t}^{t+\delta
}h(X_{s}^{\gamma_{t}},Y^{\gamma_{t}}(s),Z^{\gamma_{t}}(s))ds+\int%
_{t}^{t+\delta}Z^{\gamma_{t}}(s)dW(s).
\]
By Lemma $3.1$ and Theorem $3.2$ of Pardoux and Peng \cite{Pardoux.E 2} and
Theorem $4.4$ of Peng and Wang \cite{Peng S 3}, we deduce that%
\[%
\begin{array}
[c]{rl}
& u^{n}(\gamma_{t,t+\delta})-u^{n}(X_{t+\delta}^{\gamma_{t}})\\
= & \int_{t}^{t+\delta}D_{s}u^{n}({\gamma_{t,s}})ds-\int_{t}^{t+\delta}%
D_{s}u^{n}(X_{s}^{\gamma_{t}})ds-\int_{t}^{t+\delta}D_{x}u^{n}(X_{s}%
^{\gamma_{t}})dX^{\gamma_{t}}(s)-\frac{1}{2}\int_{t}^{t+\delta}D_{xx}%
u^{n}(X_{s}^{\gamma_{t}})d\langle X^{\gamma_{t}}\rangle(s),
\end{array}
\]
Thus by the dominated convergence theorem, we have%
\[%
\begin{array}
[c]{rl}
& u(\gamma_{t,t+\delta})-u(\gamma_{t})\\
= & -\int_{t}^{t+\delta}D_{x}u(X_{s}^{\gamma_{t}})dX^{\gamma_{t}}(s)-\frac
{1}{2}\int_{t}^{t+\delta}D_{xx}u(X_{s}^{\gamma_{t}})d\langle X^{\gamma_{t}%
}\rangle(s)\\
& +\int_{t}^{t+\delta}h(X_{s}^{\gamma_{t}},Y^{\gamma_{t}}(s),Z^{\gamma_{t}%
}(s))ds+\int_{t}^{t+\delta}Z^{\gamma_{t}}(s)dW(s)+\lim_{n\rightarrow\infty
}C^{n},
\end{array}
\]
where
\[
C^{n}=\int_{t}^{t+\delta}D_{s}u^{n}({\gamma_{t,s}})ds-\int_{t}^{t+\delta}%
D_{s}u^{n}(X_{s}^{\gamma_{t}})ds.
\]
Note that $u^{n}(\gamma_{t})\in\mathbb{C}_{l,lip}^{0,2}(\Lambda)T.$By Lemma
$3.1$ and $3.3$, we get
\[
\mid D_{s}u^{n}(\gamma_{t,s})-D_{s}u^{n}(X_{s}^{\gamma_{t}})\mid\leq
c\parallel\gamma_{t,s}-X_{s}^{\gamma_{t}}\parallel,\text{ \ \ }a.e.s\in\lbrack
t,T]\text{\ \ }P-a.s
\]
for some constant c depending on $C,T,\gamma_{t}$ and $k$. Hence
\[
\mid C^{n}\mid\leq c\delta\sup_{s\in{[t,t+\delta]}}\mid X^{\gamma_{t}%
}(s)-\gamma_{t}(t)\mid.\text{ \ }P-a.s\text{\ }%
\]
Taking expectation on both sides of $(4.2)$, we have
\[
\lim_{\delta\rightarrow0}\frac{u(\gamma_{t,t+\delta})-u(\gamma_{t})}{\delta
}=-\mathcal{L}u(\gamma_{t})+h(\gamma_{t},u(\gamma_{t}),D_{x}u(\gamma
_{t})\sigma(\gamma_{t})).
\]
Thus $u(\gamma_{t})$ belongs to $\mathbb{C}^{1,2}(\Lambda)$ and satisfies the
equation $(4.3)$. \quad$\Box$

%Under Assumption (2.1), (2.2) and $\sigma(\gamma_t,u(\gamma_t),v(\gamma_t))=\sigma(\gamma_t,u(\gamma_t))$  which doesn't dependent on $v(\gamma_t)$, we shall show that the value function $u(\gamma_t)=Y^{\gamma_t}(t))$ is actually a viscosity solution of the next path-dependent equation.
%\begin{align*}
%&D_tu(\gamma_t)+\mathcal{L}u(\gamma_t)-h(\gamma_t,u(\gamma_t),D_xu(\gamma_t)\sigma(\gamma_t,u(\gamma_t)))=0,\\\tag{5.4}
%\end{align*}
%%%%%%%%%%%%%%%%%%%%%%%参考文献
\renewcommand{\refname}

\end{document}